\documentclass[10pt, a4paper]{amsart}
\usepackage{amsmath, amsfonts, amssymb, mathtools, url}
\usepackage[pdftex]{graphicx}
\usepackage{enumitem}
\usepackage[usenames, dvipsnames]{color}
\usepackage[colorlinks=true,
        raiselinks=true,
        linkcolor=MidnightBlue,
        urlcolor=PineGreen,
        citecolor=Brown,
        pdfauthor={Who knows?},
        pdftitle={Lp bounds under quantized control},
        pdfkeywords={stochastic stability, quantization},
        pdfsubject={Technical Report},
        plainpages=false]{hyperref}

\usepackage{txfonts}
\mathtoolsset{mathic=true}
\allowdisplaybreaks
\newtheorem{theorem}{Theorem}

\newtheorem{proposition}[theorem]{Proposition}
\theoremstyle{definition}

\theoremstyle{remark}
\newtheorem{assumption}[theorem]{Assumption}

\DeclareMathOperator{\rank}{rank}

\DeclareMathOperator{\sat}{sat}

\newcommand{\R}{\mathbb{R}}

\newcommand{\Nz}{\mathbb{N}_{0}}

\newcommand{\Let}{\coloneqq}

\newcommand{\EE}{\mathsf E}
\newcommand{\PP}{\mathsf P}
\newcommand{\sigalg}{\mathfrak F}

\newcommand{\transp}{^{\mathsf T}}
\newcommand{\inverse}{^{-1}}

\newcommand{\lra}{\longrightarrow}

\renewcommand{\ge}{\geqslant}
\renewcommand{\le}{\leqslant}

\newcommand{\abs}[1]{\left\lvert{#1}\right\rvert}
\newcommand{\norm}[1]{\left\lVert{#1}\right\rVert}

\newcommand{\AssumptionEnd}{\hspace{\stretch{1}}{$\diamondsuit$}}

\newcommand{\st}{x}

\newcommand{\inp}{u}

\newcommand{\stnoise}{w}

\newcommand{\A}{A}
\newcommand{\B}{B}

\newcommand{\quan}{\mathfrak{q}}
\newcommand{\reachindex}{\kappa}
\newcommand{\binpivots}{Q}
\newcommand{\filtration}{\mathfrak{F}}
\newcommand{\proj}{\Pi}

\newcommand{\reachab}{\mathcal{R}}

\title[On mean-square boundedness under quantized observations]{On mean-square boundedness of stochastic linear systems with quantized observations}
\thanks{This research was partially supported by the Swiss National Science Foundation under grant 200021-122072 and by the European Commission under the project Feednetback FP7-ICT-223866 (www.feednetback.eu).}
\author[D.\ Chatterjee]{Debasish Chatterjee}
\author[P.\ Hokayem]{Peter Hokayem}
\author[F.\ Ramponi]{Federico Ramponi}
\author[J.\ Lygeros]{John Lygeros}
\address{D. Chatterjee, P. Hokayem, F. Ramponi, J. Lygeros are with the Automatic Control Laboratory, ETH Z\"urich, Physikstrasse 3, 8092 Z\"urich, Switzerland.}
\email{\texttt{\{chatterjee,hokayem,ramponif,lygeros\}@control.ee.ethz.ch}}

\begin{document}

	\maketitle

	\vspace*{-0.5in}

	\begin{abstract}
		We propose a procedure to design a state-quantizer with \emph{finite} alphabet for a marginally stable stochastic linear system evolving in $\R^d$, and a bounded policy based on the resulting quantized state measurements to ensure bounded second moment in closed-loop.
	\end{abstract}

	\section{Introduction and Result}
		Consider the linear control system
		\begin{equation}
		\label{e:sys}
			\st_{t+1}	= \A \st_t + \B \inp_t + \stnoise_t,\quad \st_0 \text{ given},\quad t=0, 1, \ldots,\tag{$\ast$}
		\end{equation}
		where the state $\st_t\in\R^d$ and the control $\inp_t\in\R^m$, $(\stnoise_t)_{t\in\Nz}$ is a mean-zero sequence of noise vectors, and $\A$ and $\B$ being matrices of appropriate dimensions. It is assumed that instead of perfect measurements of the state, quantized state measurements are available by means of a quantizer $\quan:\R^d\lra \binpivots$, with $\binpivots\subset\R^d$ being a set of vectors in $\R^d$ called alphabets/bins.

		Our objective is to construct a quantizer with finite alphabet and a corresponding control policy such that the magnitude of the control is \emph{uniformly bounded}, i.e., for some $U_{\max} > 0$ we have $\norm{\inp_t} \le U_{\max}$ for all $t$, the number of alphabets $\binpivots$ is \emph{finite}, and the states of \eqref{e:sys} are \emph{mean-square bounded} in closed-loop.

		Stabilization with quantized state measurements has a rich history, see e.g., \cite{ref:Del-90, ref:NaiFagZamEva-07, ref:BroLib-00, ref:EliaMitter-2002TAC, ref:TatSahMit-04, ref:Yuk-10} and the references therein. While most of the literature investigates quantization techniques for stabilization under communication constraints, especially of systems with eigenvalues outside the closed unit disc, our result is directed towards ``maximally coarse'' quantization---with finite alphabet of Lyapunov stable systems; communication constraints are not addressed in this work. The authors are not aware of any prior work dealing with stabilization with finite alphabet in the context of Lyapunov stable systems. Observe that unlike deterministic systems, local stabilization of stochastic systems with unbounded noise, at least one eigenvalue with magnitude greater than $1$, and bounded inputs is impossible.

		\begin{assumption}
		\label{a:basic}
			\mbox{}
			\begin{itemize}[label=$\circ$, leftmargin=*]
				\item The matrix $\A$ is Lyapunov stable---the eigenvalues of $\A$ have magnitude at most $1$, and those on the unit circle have equal geometric and algebraic multiplicities.
				\item The pair $(\A, \B)$ is reachable in $\reachindex$ steps, i.e., $\rank\begin{pmatrix}\B & \A\B & \cdots & \A^{\reachindex-1}\B\end{pmatrix} = d$.
				\item $(\stnoise_t)_{t\in\Nz}$ is a mean-zero sequence of mutually independent noise vectors satisfying \(C_4 \Let \sup_{t\in\Nz}\EE\bigl[\norm{w_t}^4\bigr] < \infty\).
				\item $\norm{u_t} \le U_{\max}$ for all $t\in\Nz$.\AssumptionEnd
			\end{itemize}
		\end{assumption}

		The policy that we construct below belongs to the class of $\reachindex$-history-dependent policies, where the history is that of the quantized states. We refer the reader to our earlier article \cite{ref:RamChaMilHokLyg-10} for the basic setup, various definitions, and in particular to \cite[\S3.4]{ref:RamChaMilHokLyg-10} for the details about a change of basis in $\R^d$ that shows that it is sufficient to consider $A$ orthogonal. We let $\reachab_k(\A, M) \Let \begin{pmatrix}\A^{k-1} M & \cdots & \A M & M\end{pmatrix}$ for a matrix $M$ of appropriate dimension, $M^+ \Let M\transp(M M\transp)\inverse $ denote the Moore-Penrose pseudoinverse of $M$ in case the latter has full row rank, and $\sigma_{\min}(M), \sigma_{\max}(M)$ denote the minimal and maximal singular values of $M$, resp. $I$ denotes the $d\times d$ identity matrix. For a vector $v\in\R^d$, let $\proj_v(\cdot) \Let \bigl\langle{\cdot},\tfrac{v}{\norm{v}}\bigr\rangle\tfrac{v}{\norm{v}}$ and $\proj_v^\perp(\cdot) \Let I - \proj_{v}(\cdot)$ denote the projections onto the span of $v$ and its orthogonal complement, respectively. For $r > 0$ the radial $r$-saturation function $\sat_r:\R^d\lra\R^d$ defined as $\sat_r(y) \Let \min\{r, \norm{y}\} \frac{y}{\norm{y}}$, and let $B_r\subset\R^d$ denote the open $r$ ball centered at $0$ with $\partial B_r$ being its boundary. We have the following theorem:

		\begin{theorem}
		\label{t:main}
			Consider the system \eqref{e:sys}, and suppose that Assumption \ref{a:basic} holds. Assume that the quantizer is such that there exists a constant \( r \) satisfying:
			\begin{enumerate}[label={\rm \alph*)}, leftmargin=2em, align=right, widest=2]
				\item \label{t:main:maxangle} \(\displaystyle{r > \frac{\sqrt{\reachindex}\,\sigma_{\max}(\reachab_\reachindex(\A, I)) \sqrt[4]{C_4}}{\cos(\varphi) - \sin(\varphi)}}\), where \(\varphi\in[0, \pi/4[\) is the maximal angle between \(z\) and \(\quan(z)\), \(z\not\in B_r\), and
				\item \label{t:main:radial} \(\quan(z) = \quan(\sat_r(z))\in\partial B_r\) for every \(z\not\in B_r\).
			\end{enumerate}
			Finally assume that that \(U_{\max} \ge r/\sigma_{\min}(\reachab_\reachindex(\A, \B))\). Then successive $\reachindex$-step applications of the control policy
			\[
				\begin{pmatrix} \inp_{\reachindex t} \\ \vdots \\ \inp_{\reachindex(t + 1)-1}\end{pmatrix} \Let -\reachab_\reachindex(\A, \B)^+ \A^\reachindex \quan(\st_{\reachindex t}),\quad t\in\Nz,
			\]
			ensures that \(\sup_{t\in\Nz} \EE_{\st_0}\bigl[\norm{\st_t}^2\bigr] < \infty\).
		\end{theorem}

\begin{figure}
		\begin{center}
			\includegraphics[scale=0.9]{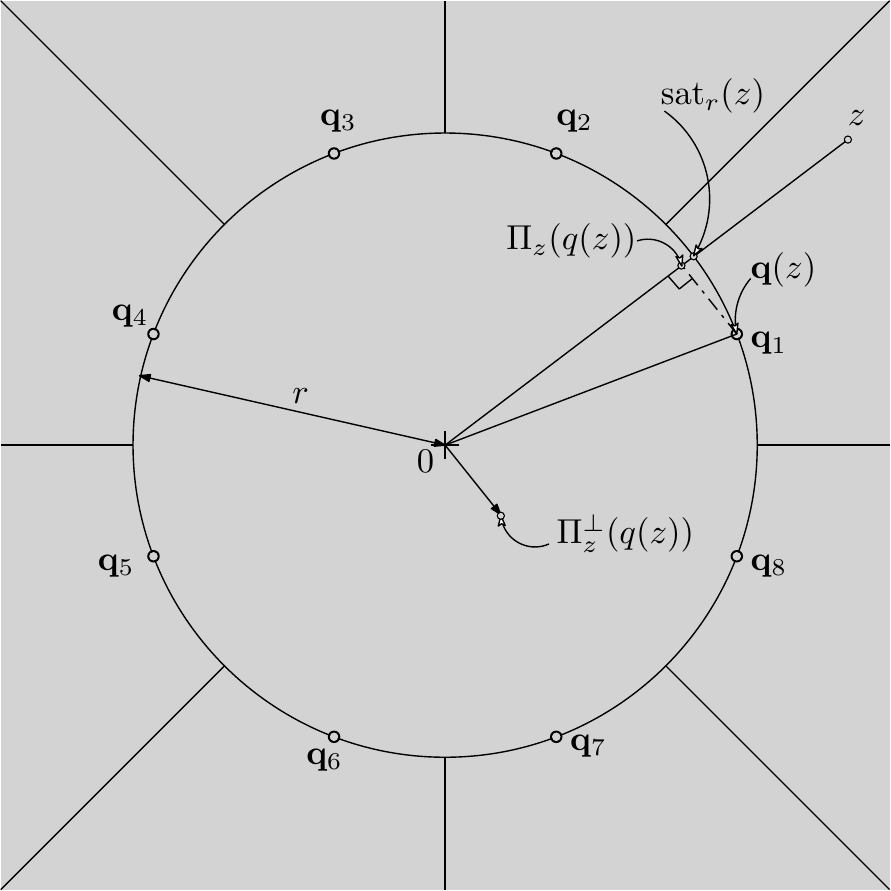}
			\caption{Pictorial depiction of the proposed quantization scheme in $\R^2$, with $\{\mathbf{q}_0 = 0, \mathbf{q}_1, \ldots, \mathbf{q}_8\}$ being the set of bins. The various projections are computed for a generic state $z$ outside the $r$-ball centered at the origin.}
			\label{fig:quan}
		\end{center}
		\end{figure}

		Observe that Theorem \ref{t:main} outlines a procedure for constructing a quantizer with finitely many bins, an example of which on $\R^2$ is depicted in Figure \ref{fig:quan}. We see from the hypotheses of Theorem \ref{t:main} that the quantizer has no large gap between the bins on the $r$-sphere, and is ``radial''; the quantization rule for states inside $B_r$ does not matter insofar as mean-square boundedness of the states is concerned. As a consequence of the control policy in Theorem \ref{t:main}, the control alphabet is also finite with $\reachindex \abs{\binpivots}$ elements. Moreover, note that as $\varphi\searrow 0$, i.e., as the ``density'' of the bins on the $r$-sphere increases, we recover the policy proposed in \cite{ref:RamChaMilHokLyg-10}, and in particular, the lower bound on $U_{\max}$ in \cite{ref:RamChaMilHokLyg-10}.

	\section{Proof of Theorem \ref{t:main}}
		We assume that the random variables $\stnoise_t$ are defined on some probability space $(\Omega, \sigalg, \PP)$. Herereafter $\EE^{\sigalg'}[\cdot]$ denotes conditional expectation for a $\sigma$-algebra $\sigalg'\subset\sigalg$.

		We need the following immediate consequence of \cite[Theorem 1]{ref:PemRos-99}.
		\begin{proposition}
		\label{p:PemRos-99}
			\textsl{
				Let $(\xi_t)_{t\in\Nz}$ be a sequence of nonnegative random variables on some probability space $(\Omega, \sigalg, \PP)$, and let $(\sigalg_t)_{t\in\Nz}$ be any filtration to which $(\xi_t)_{t\in\Nz}$ is adapted. Suppose that there exist constants $b > 0$, and $J, M < \infty$, such that $\xi_0\le J$, and for all $t$:
				\begin{gather*}
					\EE^{\sigalg_t}[\xi_{t+1} - \xi_t] \le -b\quad \text{on the event }\{\xi_t > J\},\quad\text{and}\\
					\EE\bigl[\abs{\xi_{t+1} - \xi_t}^4\big|\xi_0,\ldots, \xi_t\bigr] \le M.
				\end{gather*}
				Then there exists a constant $\gamma = \gamma(b, J, M) > 0$ such that $\displaystyle{\sup_{t\in\Nz}\EE\bigl[\xi_t^{2}\bigr] \le \gamma}$.
			}
		\end{proposition}

		\emph{Proof of Theorem \ref{t:main}: }
			Let $\filtration_t$ be the $\sigma$-algebra generated by $\{\st_s\mid s = 0, \ldots, t\}$. Since $\quan$ is a measurable map, it is clear that $(\quan(x_t))_{t\in\Nz}$ is $(\filtration_t)_{t\in\Nz}$-adapted.

			We have, for $t\in\Nz$, on $\{\norm{\st_{\reachindex t}} > r\}$,
			\begin{align*}
				\EE^{\sigalg_{\reachindex t}}\bigl[\norm{\st_{\reachindex(t+1)}} - \norm{\st_{\reachindex t}}\bigr] = \EE^{\sigalg_{\reachindex t}}\bigl[\norm{\A^{\reachindex} \st_{\reachindex t} + \reachab_\reachindex(\A, \B)\bar \inp_{\reachindex t} + \bar\stnoise_{\reachindex t}} - \norm{\st_{\reachindex t}}\bigr],
			\end{align*}
			where $\bar\inp_{\reachindex t} \Let \begin{pmatrix} \inp_{\reachindex t}\\ \vdots\\ \inp_{\reachindex(t + 1)-1}\transp\end{pmatrix}\in\R^{\reachindex m}$, and $\bar \stnoise_{\reachindex t} \Let \reachab_\reachindex(\A, I)\begin{pmatrix}\stnoise_{\reachindex t}\\ \vdots\\ \stnoise_{\reachindex(t+1)-1}\end{pmatrix}\in\R^{\reachindex d}$ is zero-mean noise. To wit, we have
			\begin{align*}
				\EE^{\sigalg_{\reachindex t}} & \bigl[\norm{\A^{\reachindex} \st_{\reachindex t} + \reachab_\reachindex(\A, \B)\bar \inp_{\reachindex t} + \bar\stnoise_{\reachindex t}} - \norm{\st_{\reachindex t}}\bigr]\\
				& \le \EE^{\sigalg_{\reachindex t}}\bigl[\norm{\A^{\reachindex} \st_{\reachindex t} + \reachab_\reachindex(\A, \B)\bar \inp_{\reachindex t}} - \norm{\st_{\reachindex t}}\bigr] + \sqrt{\reachindex}\,\sigma_{\max}(\reachab_\reachindex(\A, I)) \sqrt[4]{C_4}.
			\end{align*}
			Selecting the controls 
			\(
				\bar\inp_{\reachindex t} = -\reachab_\reachindex(\A, \B)^+ \A^{\reachindex} \quan(\st_{\reachindex t})
			\)
			as in the theorem and using the fact that $\quan(\st_{\reachindex t}) = \proj_{\st_{\reachindex t}}(\st_{\reachindex t}) + \proj_{\st_{\reachindex t}}^\perp (\st_{\reachindex t})$, we arrive at
			\begin{align*}
				\EE^{\sigalg_{\reachindex t}} & \bigl[\norm{\A^{\reachindex} \st_{\reachindex t} + \reachab_\reachindex(\A, \B)\bar \inp_{\reachindex t}} - \norm{\st_{\reachindex t}}\bigr]\\
				& = \norm{ \A^{\reachindex} \st_{\reachindex t} - \A^{\reachindex} \quan(\st_{\reachindex t}) } - \norm{ \st_{\reachindex t} }\\
				& \le \norm{ \A^{\reachindex} \st_{\reachindex t} - \sat_r(\A^{\reachindex} \st_{\reachindex t}) } - \norm{ \st_{\reachindex t} } + \norm{ \proj_{\st_{\reachindex t}}^\perp(\A^{\reachindex}\quan(\st_{\reachindex t})) }\\
				& \qquad + \norm{ \sat_r(\A^{\reachindex} \st_{\reachindex t}) - \proj_{\st_{\reachindex t}}(\A^{\reachindex}\quan(\st_{\reachindex t})) }\\
				& = -r + \norm{\A^\reachindex\sat_r(\st_{\reachindex t}) - \proj_{\st_{\reachindex t}}(\A^\reachindex\quan(\st_{\reachindex t}))} + \norm{\proj_{\st_{\reachindex t}}^\perp(\A^\reachindex\quan(\st_{\reachindex t}))}\\
				& = -r + \norm{\A^\reachindex\sat_r(\st_{\reachindex t}) - \proj_{\st_{\reachindex t}}\bigl(\A^\reachindex\quan(\sat_r(\st_{\reachindex t}))\bigr)}\\
				& \qquad + \norm{\proj_{\st_{\reachindex t}}^\perp\bigl(\A^\reachindex\quan(\sat_r(\st_{\reachindex t}))\bigr)}\;\;\text{by hypothesis \ref{t:main:radial}}\\
				& \le -r + r(1 - \cos(\varphi)) + r \sin(\varphi)\\
				& \le -b\quad\text{for some $b > 0$ by hypothesis \ref{t:main:maxangle}}.
			\end{align*}

			Moreover, we see that for $t\in\Nz$, since $\A$ is orthogonal,
			\begin{multline*}
				\EE	\Bigl[\Bigl|\norm{\st_{\reachindex(t+1)}} - \norm{\st_{\reachindex t}}\Bigr|^4 \,\Big|\,\{\norm{\st_{\reachindex s}}\}_{s=0}^t\Bigr] = \EE\Bigl[\Bigl|\norm{\st_{\reachindex(t+1)}} - \norm{\A^\reachindex \st_{\reachindex t}}\Bigr|^4 \,\Big|\,\{\norm{\st_{\reachindex s}}\}_{s=0}^t\Bigr]\\
					= \EE\bigl[\bigl|\norm{\A^\reachindex \st_{\reachindex t} + \reachab_\reachindex(\A, \B)\bar\inp_{\reachindex t} + \reachab_\reachindex(\A, I)\bar\stnoise_{\reachindex t}} - \norm{\A^\reachindex \st_{\reachindex t}}\bigr|^4 \,\big|\,\{\norm{\st_{\reachindex s}}\}_{s=0}^t\bigr]\\
					\le \EE\bigl[\norm{ \reachab_\reachindex(\A, \B)\bar\inp_{\reachindex t} + \reachab_\reachindex(\A, I)\bar\stnoise_{\reachindex t}}^4 \,\big|\,\{\norm{\st_{\reachindex s}}\}_{s=0}^t\bigr] \le M
			\end{multline*}
			for some $M > 0$ since $\bar\inp_{\reachindex t}$ is bounded in norm and $\EE\bigl[\norm{\stnoise_t}^4\bigr] \le C_4$ for each $t$.

			We let $J \Let \max\{\norm{\st_0}, r\}$. It remains to define $\xi_t \Let \norm{\st_{\reachindex t}}$ and appeal to Proposition \ref{p:PemRos-99} with the above definition of $(\xi_t)_{t\in\Nz}$ to conclude that there exists some $\gamma = \gamma(b, J, M) > 0$ such that $\sup_{t\in\Nz} \EE\bigl[\xi_t^2\bigr] = \sup_{t\in\Nz} \EE_{\st_0}\bigl[\norm{\st_{\reachindex t}}^2\bigr] \le \gamma$. A standard argument, e.g., as in \cite[Proof of Lemma 9]{ref:RamChaMilHokLyg-10}, shows that this is enough to guarantee $\sup_{t\in\Nz} \EE_{\st_0}\bigl[\norm{\st_t}^2\bigr] \le \gamma'$ for some $\gamma' > 0$.\hfill{}$\square$
		
	\section{Simulation}

		The figure below shows the average of the square norm of the state over $1000$ runs of the system:
		\begin{displaymath}
			x_{t+1} = 
			\begin{pmatrix}
			\cos(\pi/3) & -\sin(\pi/3) \\
			\sin(\pi/3) & \cos(\pi/3) 
			\end{pmatrix}
			x_t + 
			\begin{pmatrix}
			1\\0
			\end{pmatrix}
			u_t
			+w_t,
		\end{displaymath}
		where $x_0 = \bigl(10\;\;10\bigr)\transp$, $w_t\in {\mathcal N}(0, I_2)$, and where $u_t$ is chosen respectively according to the policy proposed in this article and the one proposed in \cite{ref:RamChaMilHokLyg-10}.

		\[
			\includegraphics[scale=0.45]{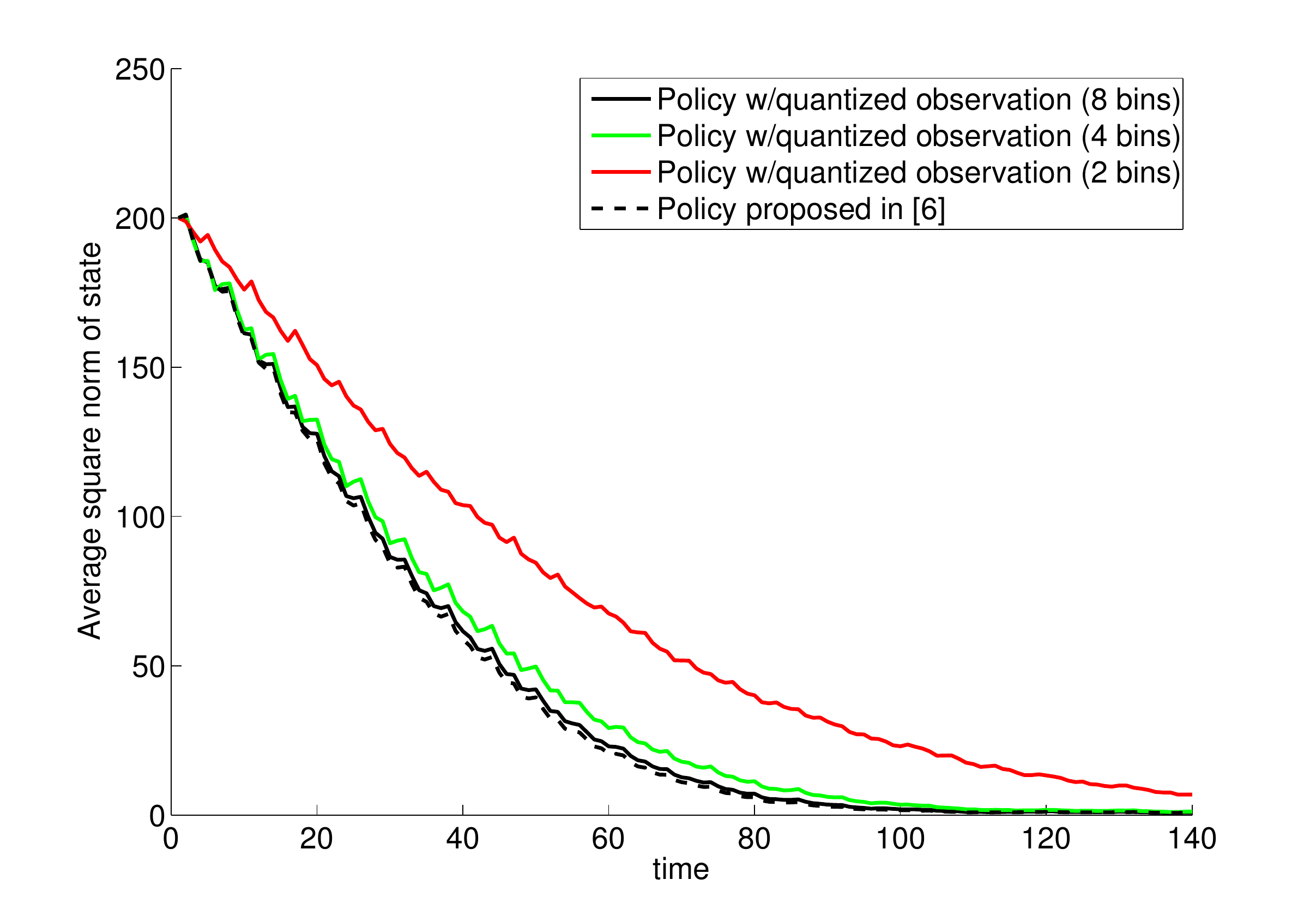}
		\]


\begin{thebibliography}{1}

\bibitem{ref:BroLib-00}
{\sc R.~W. Brockett and D.~Liberzon}, {\em Quantized feedback stabilization of
  linear systems}, IEEE Transactions on Automatic Control, 45 (2000),
  pp.~1279--1289.

\bibitem{ref:Del-90}
{\sc D.~F. Delchamps}, {\em Stabilizing a linear system with quantized state
  feedback}, IEEE Transactions on Automatic Control, 35 (1990), pp.~916--924.

\bibitem{ref:EliaMitter-2002TAC}
{\sc N.~Elia and S.~Mitter}, {\em Stabilization of linear systems with limited
  information}, Automatic Control, IEEE Transactions on, 46 (2002),
  pp.~1384--1400.

\bibitem{ref:NaiFagZamEva-07}
{\sc G.~N. Nair, F.~Fagnani, S.~Zampieri, and R.~J. Evans}, {\em Feedback
  control under data rate constraints: An overview}, Proceedings of the IEEE,
  95 (2007), pp.~108--137.

\bibitem{ref:PemRos-99}
{\sc R.~Pemantle and J.~S. Rosenthal}, {\em Moment conditions for a sequence
  with negative drift to be uniformly bounded in {$L\sp r$}}, Stochastic
  Processes and their Applications, 82 (1999), pp.~143--155.

\bibitem{ref:RamChaMilHokLyg-10}
{\sc F.~Ramponi, D.~Chatterjee, A.~Milias-Argeitis, P.~Hokayem, and
  J.~Lygeros}, {\em Attaining mean square boundedness of a marginally stable
  stochastic linear system with a bounded control input}, IEEE Transactions on
  Automatic Control, 55 (2010), pp.~2414--2418.

\bibitem{ref:TatSahMit-04}
{\sc S.~Tatikonda, A.~Sahai, and S.~Mitter}, {\em Stochastic linear control
  over a communication channel}, IEEE Transactions on Automatic Control, 49
  (2004), pp.~1549--1561.

\bibitem{ref:Yuk-10}
{\sc S.~Y{\"u}ksel}, {\em Stochastic stabilization of noisy linear systems with
  fixed-rate limited feedback}, IEEE Transactions on Automatic Control, 55
  (2010), pp.~2847--2853.

\end{thebibliography}
\end{document}